\documentclass[11pt]{article}
    \usepackage{soul}
    \usepackage{amssymb}
    \usepackage{amsfonts}
    \usepackage{eucal}
    \usepackage{amsmath}
    \newtheorem{theorem}{Theorem}[section]
    \newtheorem{proposition}[theorem]{Proposition}

    \newtheorem{lemma}[theorem]{Lemma}
    \newtheorem{remark}[theorem]{Remark}

    \newenvironment{proof}[1][Proof]{\textbf{#1.} }{\ \rule{0.5em}{0.5em}}
    \input pictex.sty

    \usepackage{cancel}
    \usepackage{xcolor}
    
    \usepackage{amsmath,color,ulem}
    \newcommand{\stkout}[1]{\ifmmode\text{\sout{\ensuremath{#1}}}\else\sout{#1}\fi}

    \DeclareMathOperator\Var{Var}

\begin{document}\title{ Intermediately Trimmed Sums of Oppenheim Expansions: a Strong Law}
 \author{Rita Giuliano\thanks{Dipartimento di
		 Matematica, Universit\`a di Pisa, Largo Bruno
		 Pontecorvo 5, I-56127 Pisa, Italy (email: \texttt {rita.giuliano@unipi.it})}
 \and~Milto Hadjikyriakou \thanks{University of Central Lancashire, 12--14 University Avenue Pyla, 7080 Larnaka, Cyprus   (email: \texttt {MHadjikyriakou@uclan.ac.uk})     
 }  
 }
  
\maketitle
\begin{abstract}
	  The work of this paper is devoted to obtaining strong laws for intermediately trimmed sums  of random variables with infinite means. Particularly, we provide conditions under which the intermediately trimmed sums of independent but not identically distributed random variables converge almost surely. Moreover, by dropping the assumption of independence we provide a corresponding convergence result for a special class of Oppenheim expansions. We highlight that the results of this paper generalize the results provided in the recent work of \cite{KS} while the convergence of intermediately trimmed sums of generalized Oppenheim expansions is studied for the first time.
\medskip

\end{abstract}


\section{Introduction} 
Consider the following framework: Let  $(\Omega, \mathcal{A}, P)$,  be the probability space  where $\Omega =[0,1]$, $\mathcal{A}$ is the $\sigma$-algebra of the Borel subsets of $[0,1]$ and $P$ is the Lebesgue measure on $[0,1]$.  Moreover, let $(B_n)_{n\geq  1}$  be a sequence of integer valued random variables  defined on $(\Omega, \mathcal{A}, P)$,   $F$ be  a  distribution function with $F(0)=0$ and $F(1)=1$,  $\varphi_n:\mathbb{N}^*\to \mathbb{R}^+$   be a sequence of  functions  and  assume that $(y_n)_{n\geq  1}$   is a sequence of nonnegative numbers  with $y_n=y_n(h_1, \dots, h_n)$ (i.e. possibly depending on the $n$ integers $h_1, \dots, h_n$) such that, for $h_1 \geq  1$ and $h_j\geq  \varphi (h_{j-1})$, $j=2, \dots, n$ we have
\begin{equation}\label{densitacondizionale}
P\big(B_{n+1}=h_{n+1}|B_{n}=h_{n}, \dots, B_{1}=h_{1}\big)=  F(\beta_n)-F(\alpha_n),
\end{equation}
where we set
$$ \delta_j(h,k, y) = \frac{ \varphi_j (h )(1+y )}{k+\varphi_j (h ) y  }, \qquad j= 1, \dots, n$$ \begin{equation*}
\alpha_n=\delta_n(h_n, h_{n+1}+1, y_n)  ,\quad \beta_n=\delta_n(h_n, h_{n+1}, y_n) .\end{equation*}

\noindent
Let $Y_n= y_n(B_1, \dots, B_n)$ and
$$R_{n}= \frac{ B_{n+1}+\varphi_n(B_n) Y_n}{\varphi_n(B_n)(1+Y_n) }= \frac{1}{\delta_n(B_n, B_{n+1}, Y_n)}.$$
In what follows the sequence $(R_n)_{n\geq  1}$ will be called an {\it Oppenheim expansion sequence  with related distribution function $F$}. We refer to the paper \cite{Gi} and to the references therein for examples of Oppenheim expansion sequences.  

\medskip

\noindent  It is known (see the seminal paper \cite{F} and also \cite{E} and \cite{K}) that when a sequence of random variables $(X_n)_{n\in \mathbb{N}}$ has infinite means, there is no normalizing sequence $(d_n)_{n\in \mathbb{N}}$ such that $\lim \dfrac{1}{d_n}\sum_{i=1}^{n}X_i = 1$ almost everywhere. On the other hand, if there is a sequence such that  $\lim \dfrac{1}{d_n}\sum_{i=1}^{n}X_i = 1$ in probability, then it might be possible to establish a strong law of large numbers after neglecting finitely many of the largest summands from the partial sum of $n$ terms (see for instance  \cite {M1}, \cite{ M2},  \cite {HMM},  \cite {M});  in such a case one says that the sum has been {\it trimmed}. The Oppenheim expansion sequences defined above are typically random variables with infinite means; nevertheless, in 
 \cite{Gi} it is proven  that, when $F$ has a density $f$, the sequence of partial sums $S_n = \sum_{k= 1}^n R_k$, normalized by $n \log n$, converges in probability to the constant 1, i.e.
$$\lim_{n \to \infty}\frac{S_n}{n \log n}=1\qquad \hbox{in probability}.$$   Moreover, in \cite{GH2021}, the same weak law has been obtained in the special case where the involved distribution functions are not continuous (see Theorem 4 and 5 therein). 
Motivated by this result, the question arises whether a strong law can be obtained for the trimmed sums of Oppenheim expansions, possibly with the same normalizing sequence  $n \log n$.
	
	\medskip
	
\noindent By using the notation introduced in    \cite{KS},  let $(X_n)_{n\geq  1}$ be a sequence of random variables and, for each $n$, let $\sigma$ be a permutation from the symmetric group $\mathcal{S}_n$ acting on ${1,2,\dots, n}$ such that $X_{\sigma(1)} \geq X_{\sigma(2)}\geq \cdots \geq X_{\sigma(n)}$. For a given sequence $(r_n)_{n\geq 1 } $ of integers with $r_n< n$ for every $n$ we set
$$S_n^{r_n}:=\sum_{k = r_n +1}^n X_{\sigma(k)}.$$ As  mentioned  above, according to the standard terminology
$(S_n^{r_n})$ is a {\it trimmed sum process}; the trimming  is called {\it light} if   $r_n = r\in \mathbb{N^*}$  for every $n$;
{\it moderate} or {\it intermediate} if $r_n\to \infty$   and $\frac{r_n}{n}\to 0$.    In \cite{GH} a  strong law of large numbers has been proved for lightly trimmed sums  of Oppenheim expansion sequences.  In the present paper we  are interested in establishing  a strong law of large numbers for intermediately trimmed sums for a special class of generalized Oppenheim expansions, and the main result of this work, Theorem \ref{teorema2}, is presented in Section 4.  

\medskip

\noindent The paper is structured as follows: Section 2 contains some preliminary results; in Section 3 we prove a strong law of large numbers for independent variables (Theorem \ref{teorema1}),   to be used for the proof of Theorem \ref{teorema2}. We point out that Theorem \ref{teorema1} is of independent interest and it can be considered as a generalization of the corresponding result presented in \cite{KS} since the convergence here is established for random variables that are not identically distributed.  The proof of the main result, Theorem \ref{teorema2}, is discussed in Section 4. It is important to highlight that  Theorem \ref{teorema2} is proven without any independence assumption since Oppenheim expansions   are not  sequences of  independent random variables  in general. Thus, Theorem \ref{teorema2} is a generalization of the classical framework of independence (as in \cite{KS}).

\medskip
\noindent Throughout the paper the notation $ \lceil x\rceil$ (resp. $\lfloor x \rfloor$) stands for the least integer greater than (resp. less than) or equal to $x$   while $a_n \sim b_n$ as $n \to \infty$,  means that $\lim_{n \to \infty}\frac{a_n}{b_n}=1$. The symbols $C$, $c$ (possibly equipped with indices) denote positive absolute constants that can take different values in different appearances,  while the notation $1_{A}(x)$ represents the indicator function for the set $A$. 

\section{Preliminaries}
	
In this section we present some results that are essential for establishing the strong laws of the following sections. We start by providing, without a proof, a well-known inequality, namely the generalized Bernstein's inequality, which can be found in \cite{H}. 

\begin{lemma} \label{lemma2}\sl    Let $Y_1, \dots, Y_n$ be independent random variables for which there exists a positive constant $M$ such that $|Y _k - EY_k|\leq M $ for every $ k =1, \dots, n$. Let $Z_n = \sum_{k=1}^n Y_k$. Then, for every $t> 0$
	$$P \Big(\max_{1\leq i \leq n}\big|Z_i - EZ_i\big|\geq t\Big)\leq 2 \exp \Big(- \frac{t^2}{2 Var (Z_n) + \frac{2}{3}Mt}\Big).$$	
\end{lemma}

\noindent  The Bernstein's inequality is the main tool for the proof of  a convergence result for the truncated random variables defined as follows. Let $(t_n)_{n \geq 1}$ be a sequence of positive numbers and let
\[
Z_n = \sum_{k=1}^n X_k 1_{\{X_k \leq t_n\}}
\] 
be the  $n-$th partial sum of the corresponding truncated   sequence. In Theorem \ref{theorem6} we  give conditions on $(t_n)_{n \geq 1}$  allowing to establish a strong law for the sequence $(Z_n)_{n \geq 1}$.

\begin{theorem} \label{theorem6}\sl Let $(X_n)_{n \geq 1} $ be a sequence  of nonnegative independent random variables and $(t_n)_{n \geq 1}$ a sequence of positive numbers. Assume that there exists $C_0  > 0$ such that
	\begin{equation}\label{assumption1}
	\sum_n \exp\Big(- C\frac{d_n^2}{n t_n^2}\Big)< + \infty 
	\end{equation}
	for all $0 <C<C_0$ where $d_n = EZ_n$.
	Then 
	$$\lim_{n \to \infty}\frac{Z_n}{d_n  }=1.$$
\end{theorem}
\begin{proof}
For simplicity, denote  $W_{k,n} = X_k 1_{\{X_k \leq t_n\}}$ for $k  =1, \dots, n$. Then, for any integer $k$ we have  
	$$\big|W_{k,n}- EW_{k,n}\big|\leq \big|W_{k,n}\big|+ \big|EW_{k,n}\big|\leq 2 t_n,$$
	$$d_n = EZ_n= \sum_{k=1}^n EW_{k,n}\leq \sum_{k=1}^n  t_n = n t_n,$$
	and
	$$Var (W_{k,n}) = E \big[(W_{k,n}- EW_{k,n})^2\big]\leq 2E\big[W^2_{k,n}+E^2W_{k,n}\big]\leq 4 t_n^2.$$
	Hence, by independence
	$$Var (Z_n) = \sum_{k=1}^n Var (W_{k,n}) \leq 4 nt_n^2 .$$
	Apply the Bernstein inequality to $W_{k,n}-EW_{k,n}$ ($ k = 1, \dots, n$) to obtain, for every $\varepsilon >0$,
	$$P\big(|Z_n -EZ_n|\geq \varepsilon EZ_n \big) \leq 2 \exp \Big(- \frac{\varepsilon^2 d_n ^2}{8 nt_n^2 + \frac{4}{3}\varepsilon t_nd_n}\Big)\leq\exp \Big(- \frac{3 \varepsilon^2}{24 + 4 \varepsilon}\cdot\frac{  d_n ^2}{  nt_n^2  }\Big).$$
	The statement follows from the Borel--Cantelli lemma, due to assumption \eqref{assumption1}. 
\end{proof}

\bigskip
\noindent
We introduce the following notation: for $n \in \mathbb{N}^*$ and $ k = 1, \dots, n$  denote 
$$a_{k,n}= P(X_k > t_n) , \qquad b_{k,n}= P(X_k \geq t_n), \qquad A_n = \sum_{k=1}^n  a_{k,n}  , \qquad B_n =\sum_{k=1}^n  b_{k,n} .$$

\medskip 

\noindent The lemma that follows can easily be obtained by employing again the generalized Bernstein's inequality.

\begin{lemma} \label{lemma1}\sl Let $(X_n)_{n \geq 1} $ be a sequence  of nonnegative independent random variables and $(t_n)_{n \geq 1}$ a sequence of positive numbers. Assume that there exists $c_0  > 0$ such that
	\begin{equation}\label{assumption2}
	\sum_n \exp\Big(- c A_n\Big)< + \infty 
	\end{equation}for all $0 <c<c_0$.
	Then, for every sufficiently small $\varepsilon >0$,
	\begin{equation}\label{statement1}
	P\left(\left|\sum_{k=1}^n 1_{\{X_k > t_n\}}- A_n  \right|\geq \varepsilon A_n  \hbox{ i.o.}\right)=0
	\end{equation}
	and
	\begin{equation}\label{statement2}
	P\left(\left|\sum_{k=1}^n 1_{\{X_k \geq t_n\}}-B_n   \right|\geq \varepsilon B_n \hbox{ i.o.}\right)=0.
	\end{equation}
\end{lemma}
\begin{proof}
	Observe that $\left|1_{\{X_k > t_n\}}-a_{k,n}  \right|= \left|1_{\{X_k > t_n\}}-  P(X_k > t_n) \right|\leq 2$ and $\Var \left(1_{\{X_k > t_n\}}-a_{k,n} \right)= a_{k,n} -a^2_{k,n} $.  By independence 
	$$\Var \left( \sum_{k=1}^n \left(1_{\{X_k > t_n\}}-a_{k,n}\right)\right) =\sum_{k=1}^n (a_{k,n} -a^2_{k,n}).$$
	Applying the Bernstein inequality we get
	\begin{align*}&
	P\left(\left|\sum_{k=1}^n 1_{\{X_k > t_n\}}- A_n  \right|\geq \varepsilon  A_n  \right)\leq 2 \exp \left(- \frac{\varepsilon^2   A_n ^2}{2 \sum_{k=1}^n (a_{k,n} -a^2_{k,n}) + \frac{4}{3}\varepsilon \sum_{k=1}^n a_{k,n} }\right)\\&\leq  \exp\left(- \frac{3 \varepsilon^2}{6 + 4 \varepsilon} A_n\right),
	\end{align*}
	which  yields that $\displaystyle\sum_{n} \exp\left(- \frac{3 \varepsilon^2}{6 + 4 \varepsilon} A_n\right)$ is finite if $\varepsilon$ is sufficiently small, by assumption \eqref{assumption2}. Thus, statement \eqref{statement1} follows from the Borel--Cantelli lemma.
	The argument is identical for statement \eqref{statement2}, since  $A_n\leq B_n  ,$   which implies
	$$ \sum_n\exp\Big(- c B_n \Big)   \leq \sum_n \exp\Big(- c A_n\Big).$$
\end{proof}

\section{A strong law for independent random variables}

\noindent
Recall the notation introduced for the needs of Lemma \ref{lemma1}, i.e. 
$$A_n = \sum_{k=1}^n P(X_k > t_n), \qquad B_n =\sum_{k=1}^n P(X_k \geq t_n) .$$
Observe that  \eqref{statement1} and \eqref{statement2} imply respectively  that, for every $\varepsilon >0$, a.e.
\begin{equation}\label{statement3}
(1-\varepsilon)A_n\leq\# \{i \leq n: X_i> t_n\}  \leq  (1+\varepsilon) A_n\quad \hbox{eventually}
\end{equation}
and
\begin{equation}\label{statement4}
(1-\varepsilon) B_n\leq \# \{i \leq n: X_i\geq  t_n\} \leq (1+\varepsilon)  B_n \quad \hbox{eventually}.
\end{equation} 
Moreover, \eqref{statement3} and \eqref{statement4} yield that, a.e.\begin{align}&\label{statement5}\nonumber
\lfloor  (1-\varepsilon) B_n\rfloor - \lceil (1+\varepsilon) A_n\rceil\leq  (1-\varepsilon)B_n- (1+\varepsilon)  A_n  \leq    \# \{i \leq n: X_i\geq  t_n\}-\# \{i \leq n: X_i>  t_n\} \\&= \# \{i \leq n: X_i= t_n\}\quad \hbox{eventually} .
\end{align}

\medskip

\noindent Before stating and proving the main theorem of this section we prove two lemmas that are essential tools for  obtaining the convergence result  we are interested in.

\begin{lemma} \label{lemma3}\sl For every pair of integers $r\geq 1$ and $n \geq 1$ we have
	$$ Z_n-S_n^{r }  =  \sum_{k=1}^{r} X_{\sigma(k)}-\sum_{k=1}^n X_k 1_{\{X_k > t_n\}}.$$
\end{lemma}
\begin{proof}
Without loss of generality we may assume that the sequence $(X_n)_{n\in\mathbb{N}}$ is already ordered in decreasing order, i.e. that $ \sigma(k)=k$: if not, just put $\tilde{X}_k = X_{\sigma(k)}$ and notice that 
$$\sum_{k=1}^n X_k 1_{\{X_k > t_n\}}= \sum_{k=1}^n \tilde X_k 1_{\{\tilde X_k > t_n\}} \mbox{ and } Z_n =  \sum_{k=1}^n \tilde X_k 1_{\{\tilde X_k \leq  t_n\}}.$$ 
Thus, we write
\begin{align*}& S_n^{r }= \sum_{k= r  +1}^n X_k= \sum_{k=  1}^n X_k- \sum_{k=  1}^{r } X_k  = \sum_{k=  1}^nX_k 1_{\{ X_k \leq  t_n\}}+\sum_{k=  1}^nX_k 1_{\{ X_k >  t_n\}}- \sum_{k=1}^{r }X_k\\&= Z_n+\sum_{k=  1}^nX_k 1_{\{ X_k >  t_n\}}- \sum_{k=1}^{r}X_k. 
\end{align*}
\end{proof}

\begin{lemma}\sl \label{lemma5}  For every  pair of integers $r\geq 1$ and $n \geq 1$  such that $r \geq  \# \{k \leq n:X_k > t_n \}$ and for every $\varepsilon > 0$ we have 
	$$Z_n- S_n^{r } \leq \big(r-(1-\varepsilon) A_n   \big)t_n.$$	
\end{lemma}
\begin{proof}
Assume again that $\sigma(k) =k$; by Lemma \ref{lemma3}, it is sufficient to prove that
\begin{align*}&
\sum_{k=1}^{r } X_{k}-\sum_{k=1}^n X_k 1_{\{X_k > t_n\}}\leq  \big(r-(1-\varepsilon) A_n   \big)t_n. 
\end{align*}
We have
\begin{align*}&
\sum_{k=1}^{r } X_{k}-\sum_{k=1}^n X_k 1_{\{X_k > t_n\}}= \sum_{k=1}^{r } X_{k}1_{\{X_k \leq t_n\}} + \sum_{k=1}^{r } X_{k}1_{\{X_k > t_n\}}-
\sum_{k=1}^n X_k 1_{\{X_k > t_n\}}=\\
&=\sum_{k=1}^{r } X_{k}1_{\{X_k \leq t_n\}}-\sum_{k=r+1}^n X_k 1_{\{X_k > t_n\}} =\sum_{k=1}^{r } X_{k}1_{\{X_k \leq t_n\}}, 
\end{align*}
since, if $k > r ,$ $X_k$ cannot be  $>t_n$ (since all the summands $>t_n$  have been trimmed before $r$, due to the assumption in the present lemma). Continuing and denoting  
$\# \{k \leq n:X_k > t_n \}= \ell_n$, we find
$$\sum_{k=1}^{r } X_{k}1_{\{X_k \leq t_n\}}= \sum_{k=1}^{\ell_n   } X_{k}1_{\{X_k \leq t_n\}}+ \sum_{k=\ell_n+1}^{r } X_{k}1_{\{X_k \leq t_n\}}= \sum_{k=\ell_n+1}^{r } X_{k}1_{\{X_k \leq t_n\}},$$
since, if $k \leq \ell_n$, then  $X_k$ is  still  $>t_n$ (by the very definition of  $\ell_n$). Then,
$$\sum_{k=\ell_n+1}^{r } X_{k}1_{\{X_k \leq t_n\}} \leq t_n\sum_{k=\ell_n+1}^{r}  1=(r -\ell_n)t_n\leq \big(r-(1-\varepsilon) A_n   \big) t_n,  $$
by the left inequality in \eqref{statement3}.
\end{proof}

\medskip 

\noindent Next, we state and prove the main result of this section which is a strong law for trimmed sums of independent random variables. As it has already been mentioned in the introduction, since the convergence result  obtained below is proven for random variables that are not identically distributed, it can be considered as a generalization of the corresponding result obtained in \cite{KS}.

\begin{theorem}\label{teorema1}\sl Let $(X_n)_{n \geq 1} $ be a sequence  of nonnegative independent random variables and $(r_n)_{n \geq 1}$ a sequence of natural  numbers tending to $\infty$ and with $r_n = o(n)$. Let $(t_n)_{n \geq 1}$ be a sequence of positive real numbers for which \eqref{assumption1} and \eqref{assumption2} hold. Assume that there exists $\varepsilon_0 > 0$ such that 
	\begin{equation}\label{assumption3}
	r_n \geq \lceil(1+\varepsilon_0)A_n \rceil 
	\end{equation}
	for  $n$ large enough. Moreover, assume that
	\begin{equation}\label{assumption4} \limsup_{n \to \infty}\frac{ A_nt_n}{d_n}=C< \infty \qquad \hbox{and}\qquad \lim_{n \to \infty}\frac{(r_n-A_n)t_n}{d_n}=0.
	\end{equation}   
	Then   $$\lim_{n \to \infty}\frac{ S_n^{r_n}}{d_n}=1.$$  
\end{theorem} 
\begin{proof}
	Notice that $Z_n\geq  S_n^{\lceil (1+\varepsilon_0) A_n \rceil }  $.  In fact, by the right inequality in \eqref{statement3},  in $S_n^{\lceil (1+\varepsilon_0) A_n \rceil }$ there are only summands not greater     than $t_n$, so the sum is not larger than the sum of all the summands not greater than $t_n$, i.e. $Z_n$. Since the assumption \eqref{assumption3} implies that $ S_n^{r_n}\leq  S_n^{\lceil (1+\varepsilon_0) A_n \rceil }$, we get  also
	\begin{equation}
	\label{statement14} S_n^{r_n}\leq Z_n.
	\end{equation}  
	Now observe that \eqref{assumption3} holds also for every $\varepsilon < \varepsilon_0$; thus 
	by the right inequality in \eqref{statement3} we have $$r_n \geq  \lceil(1+\varepsilon) A_n \rceil  \geq \# \{k \leq n: X_k> t_n\}.$$ Thus Lemma \ref{lemma5} can be applied, yielding 
	\begin{equation}\label{statement11}
	Z_n- S_n^{r_n}  \leq (r_n-(1-\varepsilon) A_n   )t_n.
	\end{equation}
	Combining relations \eqref{statement14} and \eqref{statement11} we obtain
	$$ 0 \leq  Z_n- S_n^{r_n} \leq (r_n-(1-\varepsilon ) A_n   )t_n.$$
	Finally the assumptions in \eqref{assumption4} ensure that 
	$$ 0 \leq \limsup_{n \to \infty} \frac{Z_n- S_n^{r_n}}{d_n}\leq C\epsilon;$$
	which concludes the proof due to the arbitrariness of $\varepsilon$ and by applying Theorem \ref{theorem6}.	
\end{proof}

 \section{A strong law  for a  special class of generalized Oppenheim expansions}

In this section, we provide a strong law for the trimmed sums of a special class of generalized Oppenheim expansions. Recall that Theorem \ref{teorema1} proven in the previous section, provides conditions under which a convergence result is established  for sequences of independent random variables. However, the independence assumption is not satisfied by all Oppenheim expansions and therefore we need to establish a convergence result without this constraint; thus, first we describe the class of Oppenheim expansions we  shall deal with. 

\bigskip
\noindent
Call {\it good} a strictly increasing sequence  $\Lambda = (\lambda_j)_{j \in \mathbb{N}}$    tending to $+ \infty$  with $\lambda_j\geq 1$ for every $j\geq 1$, $\lambda_0 =0$ and such that\begin{equation}\label{as1F}
\sup_n (\lambda_{n+1}- \lambda_n)= \ell< + \infty.
\end{equation} For $u \in [1, + \infty)$ let $j_u$ be the only integer such that $\lambda_{j_u-1}\leq u < \lambda_{j_u}$  i.e. $\lambda_{j_u}$ is the minimum element in $\Lambda$ greater than $u$. 

\medskip

\noindent
  The sequence defined above is instrumental for defining a class of Oppenheim expansions with optimal properties. This class is presented and studied in depth in the   paper  in progress \cite{GH}, where the proof of the following proposition can be found together with some examples.  
  
 \begin{proposition}\sl \label{indepTn} Let $(R_n)_{n \geq 1}$ be an Oppenheim expansion sequence  with related distribution function $F$; 
assume that there exists a good sequence  $\Lambda$  such that
  for every $x \in \Lambda$ and for every $n$,  $x \phi_n(B_n)+(x -1)Y_n\phi_n(B_n)$ is an integer.
  Denote
  \begin{equation}\label{Tn}
  X_n = :\lambda_{j_{R_n}}.   
  \end{equation}
    Then, the variables $X_n$ take values in $\Lambda$, are independent and $X_n$ has discrete density   given by  $$F\Big(\frac{1}{\lambda_{s-1}}\Big)-F\Big(\frac{1}{\lambda_s}\Big),\quad s\in \mathbb{N}^*.$$
 \end{proposition}

\noindent
 Next, we define  the function
$$\phi(u) := \sum_{j=2}^{j_u-1} \frac{\lambda_j - \lambda_{j-1}}{\lambda_{j-1}}$$
which plays a crucial role in the proof of the main result of this section. In the Lemma that follows we provide lower and upper bounds for this important function. 
  
 \begin{lemma}\label{lemma4} \sl For the function $\phi$ the following inequalities hold true. 
 	  \begin{enumerate}
 		 \item[(a)]   For every $u > 0$,  $$ \phi(u) \geq \log u - \log\lambda_1-\ell.$$
 		
 		\item [(b)] Let $\varepsilon >0$ be fixed. Then,  for sufficiently large $u > 0$,
 		$$  \phi(u) \leq (1+\varepsilon)\{\log u - \log\lambda_1\}.$$
 	\end{enumerate}
 \end{lemma}
\begin{proof}
 For the first inequality we start by considering the function
 $$f(x)=\sum_{j = 2}^{j_u} \frac{1}{\lambda_{j-1}}1_{[\lambda_{j-1},\lambda_j )}(x),$$
 defined on the interval $[\lambda_1, \lambda_{{j_u}}]$. Then clearly
 $$ \int_{\lambda_1}^{\lambda_{j_u}} f(x) \, {\rm d}x \geq \int_{\lambda_1}^{\lambda_{j_u}} \frac{1}{x}\, {\rm d}x =\log \lambda_{j_u}  - \log\lambda_1\geq  \log u - \log\lambda_1.$$
 And now
 $$\phi(u) = \int_{\lambda_1}^{\lambda_{j_u}} f(x) \, {\rm d}x - \frac{\lambda_{j_u}- \lambda_{j_u-1}}{\lambda_{j_u-1}}\geq \log u - \log\lambda_1 - \ell.$$
 \medskip
 
 \noindent For the second part consider the sequences
 $$C_n= \sum_{j = 2}^{n} \frac{1}{\lambda_{j}}(\lambda_j-\lambda_{j-1} ),$$
 $$D_n= \sum_{j = 2}^{n} \frac{1}{\lambda_{j-1}} (\lambda_j-\lambda_{j-1} ).$$
 Notice that, by an argument similar as in the proof of the first inequality
 $$D_n\geq \int_{\lambda_1}^{ \lambda_n}\frac{1}{x} {\rm d}x = \log \lambda_n - \log \lambda_1 \to \infty,$$
 and 
 $$C_n\leq \int_{\lambda_1}^{ \lambda_n}\frac{1}{x} {\rm d}x = \log \lambda_n - \log \lambda_1 .$$
 Thus, by applying the Cesaro  theorem, we have
 $$\lim_{n \to \infty} \frac{C_n}{D_n}=\lim_{n \to \infty}\frac{\frac{1}{\lambda_{n}}(\lambda_n-\lambda_{n-1} )}{ \frac{1}{\lambda_{n-1}} (\lambda_n-\lambda_{n-1} )}= \lim_{n \to \infty}\frac{\lambda_{n-1}}{\lambda_{n}}=1,$$
 since
 $$\frac{\lambda_{n-1}}{\lambda_{n}}= 1 - \frac{\lambda_{n}-\lambda_{n-1}}{\lambda_{n}}$$
 and
 $$0 \leq\frac{\lambda_{n}-\lambda_{n-1}}{\lambda_{n}}\leq \frac{\ell}{\lambda_{n}}\to 0, \qquad n \to \infty. $$
 As a consequence, for sufficiently large $n$ we have
 $$D_n \leq (1+\varepsilon)C_n\leq (1+\varepsilon)(\log \lambda_n - \log \lambda_1)$$
 Hence, for $n= j(u) -1$ and sufficiently large $u$, we obtain
 $$\phi(u)\leq (1+\varepsilon)(\log \lambda_{j(u) -1} - \log \lambda_1)\leq  (1+\varepsilon)(\log u - \log \lambda_1)$$
 (observe that $j(u) \to \infty$ if $u\to \infty$).
 \end{proof}
 
 \medskip
 
 \medskip
 
\noindent In this section the notation $S_n^{r_n}$ is reserved for   $(
		R_n)_{n \geq 1}$, i.e.
		$$S_n^{r_n}:=\sum_{k = r_n +1}^n R_{\sigma(k)};$$ for $X_n $ defined as in \eqref{Tn} we denote
		
		$$\tilde S_n^{r_n}:=\sum_{k = r_n +1}^n X_{\sigma(k)}.$$  
		We recall again the notation    $$ A_n:= \sum_{k =1}^n P(X_k > t_n), \qquad  d_n: =  \sum_{k =1}^nE[X_k 1_{\{X_k \leq t_n\}}] .$$
 where $(t_n)_{n \geq 1}$ is a sequence of positive numbers. For future use in this section, we investigate which further properties  $(t_n)_{n \geq 1}$ must have in order to satisfy \eqref{assumption1} and \eqref{assumption2}.

\medskip   

\begin{lemma}\label{lemma6}\sl Let $(R_n)_{n \geq 1}$ be an Oppenheim expansion sequence  with related distribution function $F$; 
assume that there exists a good sequence  $\Lambda$  such that
  for every $x \in \Lambda$ and for every $n$,  $x \phi_n(B_n)+(x -1)Y_n\phi_n(B_n)$ is an integer; assume in addition that    two positive constants $C_1 <C_2$ exist such that 
	$$C_1 \leq \frac{F(x)}{x}\leq C_2 \qquad \forall \, x \in (0,1].$$
Let $(t_n)_{n \geq 1}$ be a nondecreasing  sequence of positive numbers such that $\displaystyle\lim_{n \to \infty} t_n = + \infty$;  then 
\begin{enumerate}
\item[(i)]for every $n \geq 1$, 
\begin{equation}\label{statement6}
  d_n \geq  C n \log t_n;
\end{equation}
\item[(ii)]assumption \eqref{assumption1} is satisfied if\begin{equation}\label{assumption5}
	\sum_n \exp \left (-C \,\frac{n \log^2 t_n}{t_n ^{2  }}\right)< + \infty,
	\end{equation}
	 for sufficiently small $C >0;$  
	\item[(iii)] assumption \eqref{assumption2} is satisfied if
	\begin{equation}\label{assumption6}\sum_n \exp \left(-  C\, \frac{n  }{t_n  }\right)< +\infty, \end{equation}
	for sufficiently small $C >0.$  
\end{enumerate}
\end{lemma}
\begin{proof} Interpreting $F(\frac{1}{0})$ as 1,  we have ultimately
	$$ d_n= - \sum_{k =1}^n \sum_{j=1}^{ j_{t_n}-1}\lambda_j\left[   F\left(\frac{1}{\lambda_j}\right)- F\left(\frac{1}{\lambda_{j-1}}\right) \right].$$
	Recalling Abel's summation formula, i.e.
	$$\sum_{j=m}^r a_j (b_{j+1}- b_j)= a_r b_{r+1}- a_mb_m -\sum_{j=m+1}^rb_j (a_j - a_{j-1}),  $$
	where we take
	$$a_j = \lambda_j,  \qquad b_j=F\left(\frac{1}{\lambda_{j-1}}\right), \qquad m =1, \qquad r =  j_{t_n}-1 , $$
	 the expression  $d_n$ becomes 
	\begin{align*}& \nonumber
	 d_n=-\sum_{k =1}^{n}\left[\lambda_{j_{t_n}-1} F\left(\frac{1}{\lambda_{j_{t_n}}}\right)- \lambda_1 - \sum_{j=2}^{ {j_{t_n}-1}} F\left(\frac{1}{\lambda_{j-1}}\right) (\lambda_j-\lambda_{j-1})\right]\\&= \nonumber n \left(\sum_{j=1}^{ {j_{t_n}-1}} F\left(\frac{1}{\lambda_{j-1}}\right)(\lambda_j-\lambda_{j-1}) -\lambda_{j_{t_n}-1} F\left(\frac{1}{\lambda_{j_{t_n}}} \right) \right)\\& \geq n\left(C_1 \sum_{j=2}^{j_{t_n}-1 }\frac{\lambda_j - \lambda_{j-1}}{\lambda_{j-1}}+\lambda_1- C_2\right)  =n\left(C_1 \phi(t_n)+\lambda_1- C_2\right)\geq  C_3 n \log t_n, 
	\end{align*} for every $n \geq 1 $ by the first part of  Lemma \ref{lemma4} (the first inequality holds since  $\lambda_{j_{t_n}} \geq\lambda_{j_{t_n}-1}$).
	Then
	\begin{align}& \nonumber
	\frac{ d_n^2}{n t_n^2} = \frac{n}{t_n ^{2  }}\left(    \sum_{j=1}^{ {j_{t_n}-1}} F\left(\frac{1}{\lambda_{j-1}}\right)(\lambda_j-\lambda_{j-1}) -\lambda_{j_{t_n}-1} F\left(\frac{1}{\lambda_{j_{t_n}}} \right)  \right)^2 \geq C_3 \,\frac{n \log^2 t_n}{t_n ^{2  }}.\end{align}
	Hence  assumption \eqref{assumption1} is satisfied  if \eqref{assumption5} holds  for sufficiently small $C >0$. 

	\medskip
	\noindent
	Concerning assumption \eqref{assumption2}, note that we can write 
	$$ A_n = \sum_{k =1}^n \sum_{j=j_ {t_n}  }^{+ \infty}\left( F\left(\frac{1}{\lambda_{j-1}}\right)-  F\left(\frac{1}{\lambda_j}\right)\right)= n F\left(\frac{1}{\lambda_{j_ {t_n} -1}}\right)\geq C_1\frac{n}{\lambda_{j_ {t_n} -1}}\geq C_1 \frac{n}{t_n},$$
	thus we obtain  assumption \eqref{assumption2} by imposing that \eqref{assumption6} holds  for sufficiently small $C >0$. 
	\end{proof}
	
	\bigskip
	\noindent The main result of this section (and of the paper) is presented below; for its proof we will employ Theorem \ref{teorema1}, in which we take $X_n =\lambda_{j_{R_n}} $. Furthermore, for the purpose of using the same Theorem, we need  a  sequence   $(t_n)_{n \geq 1} $  which satisfies assumptions \eqref{assumption1} and \eqref{assumption2}.   Based on the discussion presented in Lemma \ref{lemma6}, a suitable sequence is $t_n = n^\gamma$, with   $0< \gamma < \frac{1}{2}$; thus in the sequel  we take $t_n = n^\gamma$ and $d_n$ means
	  $$ d_n =  \sum_{k =1}^n E[X_k 1_{\{X_k \leq n^\gamma\}}].$$
 
\begin{theorem} \label{teorema2} \sl Let $(R_n)_{n \geq 1}$ be an Oppenheim expansion sequence with related   distribution function $F$.  Assume that there exists a good sequence $\Lambda$    such that, for every $x \in \Lambda$ and for every $n$,  $x \phi_n(B_n)+(x -1)Y_n\phi_n(B_n)$ is an integer.   
	
	\begin{enumerate}
    \item[ (a)] Assume that there exist two positive constants $C_1 <C_2$ with
	$$C_1 \leq \frac{F(x)}{x}\leq C_2 \qquad \forall \, x \in (0,1].$$
	Then there exists a constant $C>0$ such that, letting  $r_n = \lceil \beta n^{1-\gamma}\rceil$,  
	where $0< \gamma < \frac{1}{2}$ and $\beta >C$, we have  
	$$ \lim_{n \to \infty} \frac{S^{r_n}_n}{ d_n}= 1.$$ 
	\item [(b)] Assume that
	$$\lim_{x \to 0} \frac{F(x)}{x}=\alpha >0.$$
	Then, for the same sequence $(r_n)_{n \geq 1}$ as in part (a), we have
	$$ \lim_{n \to \infty} \frac{S^{r_n}_n}{n \log n}= \alpha\gamma  .$$
	\end{enumerate}
\end{theorem}
	\begin{proof} 
	\noindent
	 For part (a), we start by observing that for every $n$,
		$$X_n - \ell \leq R_n\leq X_n,$$ which leads to 
		\[
		\sum_{k=r_n+1}^{n}(X_{\sigma(k)}-\ell)\leq \sum_{k=r_n+1}^{n}R_{\sigma(k)}\leq \sum_{k=r_n+1}^{n}X_{\sigma(k)}.
		\]	
Thus
	$$\tilde S_n^{r_n}- \ell n\leq \tilde S_n^{r_n}- \ell n+\ell r_n  \leq  S_n^{r_n}\leq \tilde S_n^{r_n},$$
	due to \eqref{as1F}. Hence, it suffices to prove   that $$ \lim_{n \to \infty} \frac{\tilde S^{r_n}_n}{d_n }=1\qquad \hbox{\rm and} \qquad  \lim_{n \to \infty}\frac{n}{d_n  }=0. $$
	The second relation is an easy consequence of \eqref{statement6}. For the first one, note that 
	\begin{equation}\label{statement8}
	(n-1)^\gamma < \lambda_{j_ {t_{n-1} }}\leq \lambda_{j_ {t_n} }\leq C_4\lambda_{j_ {t_{n }-1 }},
	\end{equation}
  (recall that $\lim_{n \to \infty}\frac{\lambda_n}{\lambda_{n-1}}=1,$ see the proof of Lemma  \ref{lemma4})
whence  
$$ A_n =nF\left(\frac{1}{\lambda_{j_ {t_n} -1}}\right)\leq C_2\frac{n}{\lambda_{j_ {t_n} -1}} \leq C_5\frac{n}{(n-1)^\gamma }\leq C_6 n^{1-\gamma} ,$$ for sufficiently large $n$.  
So for any $\varepsilon_0>0$ we have that 
\[
(1+ \varepsilon_0)A_n\leq C_6(1+ \varepsilon_0)n^{1- \gamma}  ,
\]  ultimately.

\bigskip
	\noindent Let $C=C_6$; choose $\beta > C$ and  $\varepsilon_0 \leq \frac{\beta -  C}{  C}$; then $\beta \geq  C(1+ \varepsilon_0)  $ and
$$r_n = \lceil\beta n^{1-\gamma}\rceil\geq \beta n^{1-\gamma} \geq C(1+ \varepsilon_0)n^{1- \gamma}\geq (1+ \varepsilon_0)A_n$$
which proves that assumption \eqref{assumption3} is satisfied.  Furthermore,
$$0 \leq \frac{  A_n t_n}{d_n}\leq \frac{Cn^{1- \gamma}n^\gamma}{n \left(\sum_{j=1}^{ {j_{t_n}-1}} F\left(\frac{1}{\lambda_{j-1}}\right)(\lambda_j-\lambda_{j-1}) -\lambda_{j_{t_n}-1} F\left(\frac{1}{\lambda_{j_{t_n}}} \right) \right)} \leq \frac{ C }{C_3 \log t_n}= \frac{C_7}{\log n}\to 0, \qquad n \to \infty$$
and
	\begin{align*}&
	0 \leq\frac{(r_n -  A_n)t_n }{ d_n }\leq  \frac{\left(\lceil\beta n^{1-\gamma} \rceil- nF\left(\frac{1}{  \lambda_{j_{t_n}-1}}\right)\right)n^\gamma}{n\left(C_1 \displaystyle\sum_{j=2}^{j_{t_n}-1} \frac{\lambda_j-\lambda_{j-1}}{\lambda_{j-1}} + \lambda_1- C_2\right)}\leq\frac{1}{C_8\log n }\left(\frac{\lceil\beta n^{1-\gamma} \rceil}{n^{1-\gamma}}-C_1\frac{n^\gamma}{  \lambda_{j_{t_n}-1}}\right)\to 0, \quad n \to \infty,   
	\end{align*}
 by \eqref{statement6} and since the term in parenthesis is bounded {  (recall that $ \frac{(n-1)^{\gamma}}{C_4} \leq \lambda_{j_{t_n}-1} \leq n^{\gamma}$, see the last inequality in \eqref{statement8})}. 
	An application of Theorem \ref{teorema1} concludes the proof of part (a).

	\bigskip
	\noindent
	For the proof of  part (b), it suffices to use  part (a) and notice that
	\begin{equation}\label{statement7}
	\lim_{n \to \infty}\frac{d_n}{n \log n}= \alpha \gamma.
	\end{equation}
	In order to prove this relation, let $\epsilon > 0$; then, for sufficiently small $x$, we have
	$$(\alpha - \epsilon)x \leq F(x) \leq (\alpha + \epsilon)x.$$
	 Using these inequalities for estimating $d_n$ { and the inequalities proved in Lemma \ref{lemma4}}, in a similar way as in \eqref{statement6} we have that {  for sufficiently large $n$,}
	 \begin{align*}
	    &d_n =n \left(\sum_{j=1}^{ {j_{t_n}-1}} F\left(\frac{1}{\lambda_{j-1}}\right)(\lambda_j-\lambda_{j-1}) -\lambda_{j_{t_n}-1} F\left(\frac{1}{\lambda_{j_{t_n}}} \right) \right){ \leq n(\alpha+\epsilon)\sum_{j=2}^{j_{t_n}-1}\dfrac{\lambda_j-\lambda_{j-1}}{\lambda_{j-1}}+\lambda_1}\\
	    &\leq n(\alpha+\epsilon)\phi(t_n)+n(\alpha+\epsilon)\lambda_1  \leq n(\alpha+\epsilon)(1+\epsilon)\log t_n+n(\alpha+\epsilon)\lambda_1\\&=
  (\alpha + \epsilon)(1+\epsilon)\gamma n \log n+n(\alpha+\epsilon)\lambda_1
	    \end{align*}
	    and
	    \begin{align*}
	    &d_n =n \left(\sum_{j=1}^{ {j_{t_n}-1}} F\left(\frac{1}{\lambda_{j-1}}\right)(\lambda_j-\lambda_{j-1}) -\lambda_{j_{t_n}-1} F\left(\frac{1}{\lambda_{j_{t_n}}} \right) \right) \geq {  n((\alpha-\epsilon)\phi(t_n) +(\alpha-\epsilon)\lambda_1-(\alpha+\epsilon))}\\
	    &{\geq n((\alpha-\epsilon)\gamma\log n-(\alpha-\epsilon)\log\lambda_1-(\alpha-\epsilon)\ell+(\alpha-\epsilon)\lambda_1-(\alpha+\epsilon))}.
	   \end{align*}
	The statement \eqref{statement7} follows immediately by the arbitrariness of $\epsilon$.	
\end{proof}

\begin{remark} Observe that in  part (b) of Theorem \ref{teorema2} the normalization sequence is $n\log n$, which is the one used in \cite{Gi} and \cite{GH2021}  for obtaining a weak law. 
\end{remark}

\end{document}